\documentclass{amsart}

\usepackage{amsmath}
\usepackage{amssymb}
\usepackage{color}
\usepackage{amsthm}
\usepackage{graphicx}
\usepackage{hyperref}

\theoremstyle{plain}
\newtheorem{theorem}{Theorem}[section]
\theoremstyle{definition}
\newtheorem{definition}{Definition}[section]
\theoremstyle{plain}
\newtheorem{prop}{Proposition}[section]
\newtheorem{remark}{Remark}[section]
\theoremstyle{plain}
\newtheorem{cor}{Corollary}[section]

\definecolor{Noir}{rgb}{0,0,0} 
\definecolor{Blanc}{rgb}{1,1,1} 
\definecolor{Gray}{rgb}{0.5,0.5,0.5} 
\definecolor{Rouge}{rgb}{0.8,0.1,0.1} 
\definecolor{DBleu}{RGB}{51,51,178} 
\definecolor{LBleu}{rgb}{0.85,0.85,1} 
\definecolor{Orange}{RGB}{255,140,0}

\newcommand{\bdoc}{\begin{document}} 
\newcommand{\edoc}{\end{document}} 
 
\newcommand{\bcent}{\begin{center}} 
\newcommand{\ecent}{\end{center}} 
 
\newcommand{\benum}{\begin{enumerate}} 
\newcommand{\eenum}{\end{enumerate}} 
\newcommand{\bitem}{\begin{itemize}} 
\newcommand{\eitem}{\end{itemize}} 
 
\newcommand{\btab}{\begin{tabular}} 
\newcommand{\etab}{\end{tabular}} 
\newcommand{\beqn}{\begin{eqnarray}} 
\newcommand{\eeqn}{\end{eqnarray}} 
 
\newcommand{\bmath}{\begin{math}} 
\newcommand{\emath}{\end{math}} 
 
\newcommand{\noin}{\noindent} 
 
\providecommand{\tb}[1]{\textbf{#1}} 
\providecommand{\mb}[1]{\mathbf{#1}} 
 
\newcommand{\bsh}{\backslash} 
 
\newcommand{\ds}{\displaystyle} 
 
\newcommand{\ub}{\underbrace} 
\newcommand{\ob}{\overbrace} 
\newcommand{\matr}[4]{\left \lbrack \begin{array}{cc} #1 & #2 \\
     #3 & #4 \end{array} \right \rbrack}

\providecommand{\F}[1]{\mathbb{#1}} 
 
\newcommand{\FF}{\F F} 
\providecommand{\Fn}[1]{\FF_{#1}} 
\newcommand{\Fp}{\Fn{p}} 
\newcommand{\Fq}{\Fn{q}} 
\newcommand{\Fpm}{\Fn{p^m}} 
\newcommand{\Fpn}{\Fn{p^n}} 
\newcommand{\Fpr}{\Fn{p^r}} 
 
\newcommand{\ZZ}{\mathbb{Z}} 
\providecommand{\Zn}[1]{\ZZ_{#1}} 
\providecommand{\ZnZ}[1]{\ZZ/#1\ZZ} 
\newcommand{\Zp}{\Zn{p}} 
\newcommand{\ZpZ}{\ZnZ{p}} 
 
\newcommand{\NN}{\F{N}} 
\newcommand{\QQ}{\F{Q}} 
\newcommand{\RR}{\F{R}} 
\newcommand{\CC}{\F{C}} 
\newcommand{\QQbar}{\overline{\QQ}} 
\newcommand{\Zero}{\mathbb{00}} 
 
\newcommand{\EE}{\F E} 
\newcommand{\II}{\F I} 
\newcommand{\KK}{\F K} 
\newcommand{\MM}{\F M} 
\newcommand{\XX}{\F X} 
\newcommand{\PP}{\F P} 
\newcommand{\FA}{\F A} 
\newcommand{\LL}{\F L} 
\newcommand{\HH}{\mathbb H}  
\newcommand{\FS}{\F S} 
\newcommand{\TT}{\F T}  
\newcommand{\FSig}{\F\Sig} 
\newcommand{\FDel}{\F\Del} 
 
\providecommand{\E}[1]{\hat{\F{#1}}} 
\newcommand{\EC}{\E{C}} 
 
\newcommand{\bQ}{\mathbf{Q}} 
\newcommand{\bP}{\mathbf{P}}

\newcommand{\ind}{\mbox{ind}} 
 
\newcommand{\fx}{f(x)} 
\newcommand{\gx}{g(x)} 
 
\newcommand{\x}{^\star} 
\newcommand{\xs}{^{~\star}} 
\providecommand{\U}[1]{\left(#1\right)\x} 
\newcommand{\xt}{^\times} 
 
\newcommand{\iso}{\simeq} 
\newcommand{\plus}{\oplus} 
\newcommand{\Plus}{\bigoplus} 
\newcommand{\tensor}{\otimes} 
\newcommand{\Tensor}{\bigotimes} 
\newcommand{\inject}{\hookrightarrow} 
\newcommand{\linject}{\hookleftarrow} 
\newcommand{\surject}{\twoheadrightarrow} 
\newcommand{\tri}{\vartriangleleft} 
 
\renewcommand{\ker}{\mbox{ker}\;} 
\newcommand{\Imf}{\mbox{im}~\;} 
\newcommand{\img}{\mbox{im}\;} 
\newcommand{\Hom}{\mbox{Hom}} 
\newcommand{\Sym}{\mbox{Sym}} 
\newcommand{\End}{\mbox{End}\,} 
\newcommand{\Endz}{\mbox{End}_0} 
\newcommand{\Id}{\mbox{Id}} 
\newcommand{\rk}{\mbox{rk}\,} 
\newcommand{\Pic}{\mbox{Pic}} 
\newcommand{\Jac}{\mbox{Jac}} 
\newcommand{\ch}{\mbox{ch}\,} 
\newcommand{\td}{\mbox{td}\,} 
 
\providecommand{\Gal}[1]{\mbox{Gal}(#1)} 
\providecommand{\GAL}[2]{\mbox{Gal}(#1/#2)} 
\providecommand{\Sub}[1]{\mbox{Sub}(#1)} 
\providecommand{\Lat}[1]{\mbox{Lat}(#1)}

\newcommand{\dup}{d_\wedge} 
\newcommand{\drt}{d_>} 

\newcommand{\MF}{\mathfrak}

\newcommand{\Div}{\,\,|\,\,} 
\newcommand{\lcm}{\mbox{lcm}} 
 
\providecommand{\leg}[2]{\left(\frac{#1}{#2}\right)} 
\providecommand{\jac}[2]{\leg{#1}{#2}} 
\providecommand{\qdc}[2]{\left[\frac{#1}{#2}\right]} 
 
\newcommand{\w}{\omega} 
\newcommand{\W}{\Omega} 
\providecommand{\Cal}[1]{\mathcal{#1}} 
\newcommand{\CL}{\Cal L} 
\newcommand{\CO}{\Cal O} 
\renewcommand{\O}{\Cal O} 
\renewcommand{\o}{\Cal O} 
\newcommand{\co}{\Cal O} 
\newcommand{\ca}{\Cal A}
\newcommand{\CR}{\Cal R} 
\newcommand{\CE}{\Cal E} 
\newcommand{\CF}{\Cal F} 
\newcommand{\CQ}{\Cal Q} 
\newcommand{\CK}{\Cal K} 
\newcommand{\CM}{\Cal M} 
\newcommand{\CN}{\Cal N} 
\newcommand{\CDee}{\Cal D} 
\newcommand{\CP}{\Cal P} 
\newcommand{\CS}{\Cal S} 
\newcommand{\ClC}{\Cal C} 
\newcommand{\CJ}{\Cal{J}} 
\newcommand{\CA}{\Cal{A}} 
\newcommand{\CB}{\Cal{B}}

\newcommand{\Si}{\Sigma} 
 
\providecommand{\Ok}[1]{\CO_{#1}} 
\newcommand{\OK}{\Ok{K}}

\renewcommand{\epsilon}{\varepsilon} 
\newcommand{\ep}{\varepsilon} 
 
\providecommand{\abs}[1]{\left|#1\right|} 
\providecommand{\norm}[1]{\lVert#1\rVert} 
 
\newcommand{\di}{\partial} 
\providecommand{\ddy}[1]{\ds\frac{d}{d #1}} 
\providecommand{\didiy}[1]{\ds\frac{\di}{\di #1}} 
\newcommand{\ddx}{\ddy{x}} 
\newcommand{\didix}{\didiy{x}} 
\providecommand{\v}[1]{\vec{#1}} 
 
\newcommand{\ihat}{\hat{\infty}}

\providecommand{\set}[1]{\left\{#1\right\}} 
\providecommand{\lst}[1]{\ds\left[\,\,#1\,\,\right]}

\providecommand{\ip}[2]{\left( #1,#2\right)} 
\providecommand{\IP}[2]{\left\langle #1,#2\right\rangle} 
\providecommand{\bra}[1]{\left\langle\left. #1\right|\right.} 
\providecommand{\ket}[1]{\left.\left| #1\right.\right\rangle} 
\providecommand{\bkv}[4]{\left\langle\left.\tb{#1} #2\right|\tb{#3} #4\right\rangle} 
\providecommand{\bk}[2]{\bkv{}{#1}{}{#2}} 
\providecommand{\bkvop}[5]{\left\langle\tb{#1} #2\left| #5\right|\tb{#3} #4\right\rangle} 
\providecommand{\bkop}[3]{\bkvop{}{#1}{}{#2}{#3}} 
\providecommand{\comm}[2]{\left[#1,#2\right]} 
\providecommand{\expn}[1]{\left\langle #1\right\rangle} 
\providecommand{\gen}[1]{\expn{#1}} 
\newcommand{\Del}{\Delta} 
\newcommand{\Nab}{\nabla} 
\newcommand{\Sig}{\Sigma} 
\newcommand{\oline}{\overline}

\newcommand{\p}{\rho} 
\providecommand{\sprod}[2]{\left\langle #1,#2\right\rangle}

\newcommand{\Ehat}{\overline E} 
\newcommand{\Phihat}{\overline \Phi} 
\newcommand{\phihat}{\overline \phi}

\newcommand\varleq{\mathbin{\vcenter{\hbox{%
  \oalign{\hfil$\scriptstyle<$\hfil\cr 
          \noalign{\kern-.3ex} 
          $\scriptscriptstyle({-})$\cr}%
}}}} 
 
\renewcommand\subsetneq{\mathbin{\vcenter{\hbox{%
  \oalign{\hfil$\scriptstyle\subset$\hfil\cr 
          \noalign{\kern-.3ex} 
          $\scriptscriptstyle({-})$\cr}%
}}}} 

\author{Lisa Jeffrey}
\address{Department of Mathematics\\
  University of Toronto\\
  Bahen Centre, 40 St. George Street, Toronto, ON, CANADA ~M5S 2E4\vspace{10pt}}
\email{jeffrey@math.toronto.edu}

\author{Aidan Lindberg}
\address{Department of Mathematics \& Statistics\\
  McGill University\\
  Burnside Hall, 805 Sherbrooke Street West, Montr\'eal, QC, CANADA ~H3A 0B9\vspace{10pt}}
\email{aidan.lindberg@mail.mcgill.ca}

\author{Steven Rayan}
\address{Department of Mathematics \& Statistics\\
  University of Saskatchewan\\
  McLean Hall, 106 Wiggins Road, Saskatoon, SK, CANADA ~S7N 5E6\vspace{10pt}}
\email{rayan@math.usask.ca}

\title[Explicit Poincar\'e Duality for an $SU(2)$ Character Variety]{Explicit Poincar\'e Duality in the Cohomology Ring of the $SU(2)$ Character Variety of a Surface}
\date{\today}

\begin{document}

\maketitle

\noin\tb{Abstract.} We provide an explicit description of the Poincar\'e dual of each generator of the rational cohomology ring of the $SU(2)$ character variety for a genus $g$ surface with central extension --- equivalently, that of the moduli space of stable holomorphic bundles of rank 2 and odd degree.\\

\section{Introduction}

The goal of this note is to take a deeper look at some known methods of calculating homological duals of the generators of the rational cohomology ring of the $SU(2)$ character variety of a Riemann surface, which is related by the Riemann-Hilbert correspondence and the Narasimhan-Seshadri Theorem \cite{NS} to the moduli space of stable bundles of rank $2$ on the same Riemann surface.  After reviewing the structure of the cohomology ring, we examine methods of systematically identifying the Poincar\'e duals of each of its three families of generators. It is worth mentioning that certain well-known results are included and proven or sketched for completeness and to better motivate arguments in later parts of the paper, with due reference.

For each of the three families of generators in question, we provide a method for calculating its Poincar\'e dual. The dual of the degree two generator is presented as the zero section of a particular line bundle, whose existence and construction is well known. The duals of the degree three generators follow from a simple topological operation.  Finally, we show that the degree four generator may be represented by the Euler class of a bundle over the moduli space, and state how to represent the dual to this class.  We finish by demonstrating how to calculate the Poincar\'e duals using Euler classes in general.

\section{Generators of the Cohomology Ring}

This section is a review of known facts and properties regarding the rational cohomology ring of $\mathcal M_g(n,d)$, the moduli space of stable holomorphic rank $n$, degree $d$ vector bundles with fixed determinant and $n$ coprime to $d$ over a smooth, compact, connected genus $g$ Riemann surface $\Sigma$.  The generators of this ring were known to Newstead \cite{Newstead1}, with relations between generators having been known since the seminal paper of Atiyah and Bott \cite{AB}; see also \cite{MN,Newstead2}. Explicit knowledge of the Betti numbers of the ring dates even earlier, cf. Newstead \cite{Newstead} and Harder and Narasimhan \cite{HN}.\\

\indent For our purposes, we will use the homeomorphism furnished by Narasimhan-Seshadri and redefine $\mathcal{M}_g(n,d)$ to be the space of representations
 of $\pi_1(\Sigma)$ into $SU(n)$, up to 
conjugation.  While this space is not complex-analytically isomorphic to the moduli space of bundles, the two spaces have the same cohomology ring.
 Now we pick loops $$a_1, \dots , a_{g}, b_1, \dots, b_g$$ in $\Sigma$ which generate
 the fundamental group 
of $\Sigma$ in the usual way. This is the group\\
\begin{equation*}
\pi_1(\Sigma)=\langle a_1, \dots , a_{g}, b_1, \dots, b_g \hspace{1ex}|\, \prod_{j=1}^g [a_j,b_j]=1 \rangle.
\end{equation*}

From now on we restrict to $n =2$ and $d = 1$ and so for economy we will write$$\mathcal M_g\;:=\;\mathcal M_g(2,1).$$We then wish to look at representations $\rho\, : \pi_1(\Sigma) \to (SU(2))^{2g}$ which send each generator of $\pi_1(\Sigma)$ to an element of $SU(2)$. Denote these representations by $A_i=\rho(a_i) \, , \hspace{1ex} B_i=\rho(b_i)$. These are subject to the relation \\
\begin{equation*}
\prod_{j=1}^g [A_j,B_j]=I.
\end{equation*}
\\
If we remove a small disk $D$ from $\Sigma$, then this relation no longer holds on $\Sigma \setminus D$, and we simply have the free group on $2g$ generators. Now, fix $G=SU(2)$ and consider the map

\begin{align*}
\mu : \hspace{1ex} (SU(2))^{2g} &\to SU(2) \\
\mu : (A_1, \dots, A_g, B_1, \dots, B_g) &{\mapsto} \prod_{j=1}^g [A_j,B_j].
\end{align*}
\\
It can be shown that \\
\begin{equation} \label{defm}
M_g:=\mu^{-1}(-I) \subset (SU(2))^{2g}
\end{equation}
\\
 is a smooth submanifold of $(SU(2))^{2g}$ with real dimension $6g-3$, as $-I$ is a regular value of $\mu$. Elements of $M_g$ can be seen as representations of $\pi_1(\Sigma)$ having holonomy $-I$ around $\partial D$. We then pass to the associated vector bundle and mod out by the action of $SU(2)/\{\pm I\} =SO(3)$, thereby recovering $\mathcal M_g$. \\

\begin{prop} The above action is free.\end{prop}

\begin{proof} Let
$\phi=\prod_{j=1}^g [A_j,B_j]$. Suppose $H$ is a subgroup of $SU(2)$ larger than the center, and is the stabilizer of a point $(A,B)$ in $\phi^{-1}(c)$, where
 $c=- I$ is the generator of the 
center of $SU(2)$. Here, $(A,B)=(A_1, \dots, A_g, B_1, \dots, B_g)$ are representations of $\pi_1(\Sigma)$ into $SU(2)$. By assumption $(A,B) \in Z(H)$ which commutes with $H$. This means that $\phi(A,B)$ is in the commutator subgroup $[Z(H),Z(H)]$, which is the identity,
 but we had assumed $c$ was not the identity. This is a contradiction.
\end{proof}

The following is now immediate:

\begin{cor} The moduli space $\mathcal{M}_g \, = M_g/SO(3)$ is a smooth manifold of real dimension $6g-6$.\end{cor}

Now, if $\mathbb{U}$ is the universal bundle over $\mathcal{M}_g \times \Sigma$ and if parentheses denote the slant product$$(\hspace{1ex},\hspace{1ex}): H^N(\Sigma \times \mathcal{M}_g(n,d)) \times H_j(\Sigma) \to H^{N-j}(\mathcal{M}_g(n,d)),$$then the generators of the cohomology ring $H^*(\mathcal{M}_g(n,d))$ for $2\leq r \leq n$ are

\begin{align*}
f_r&=([\Sigma], c_r(\mathbb{U})) \\
b_r^j&=(\alpha_j, c_r(\mathbb{U})) \\
a_r&= (p, c_r(\mathbb{U})), 
\end{align*}
\\
\noin where $[\Sigma] \in H_2(\Sigma)$, $\alpha_j \in H_1(\Sigma)$ are the $2g$ cycles arising from the generators of $\pi_1(\Sigma)$, and $p$ is a point which provides a generator of $H_0(\Sigma)$. 
Here, $c_r$ is the $r$-th Chern class. These homology classes form a basis of $H_2(\Sigma)$ and $H_1(\Sigma)$ respectively.

In specializing to $\mathcal M_g=\mathcal{M}_g(2,1)$, the generators take the form 

\begin{align*}
a\in &  \hspace{1ex} H^4(\mathcal{M}_g) \\
b_j \in & \hspace{1ex}  H^3(\mathcal{M}_g) ~~(1\le j \le 2g)\\
f \in & \hspace{1ex}  H^2(\mathcal{M}_g) .
\end{align*} \\
 Here the classes $b_j$ are
associated to a choice of a basis for $H^1(\Sigma, \ZZ)$,
chosen so that $b_{k}$ and $b_{k+g}$ 
(for all $ 1 \le k \le g$) 
satisfy that their cup product $b_{k}\cup b_{k+g}$ is
a chosen generator of $H^2(\Sigma,\ZZ) $ corresponding to the orientation
of $\Sigma$, and the cup product of $b_k$ with $b_j$ 
 is $0$ whenever $j \ne k+g$.

We will now describe intersection pairings --- in other words, polynomials in these generators, evaluated
on the fundamental class $[\mathcal M_g]$.  Thaddeus  shows in \cite{Thaddeus}
(equation 25, p. 144) that the only potentially nonzero intersection pairings are those of the form
\begin{equation} \label{twentyfive}
(a^m f^n \gamma_{i_1} \gamma_{i_2} \dots \gamma_{i_p}) [\mathcal{M}_g].
\end{equation}\\
where $1\leq i_1 < i_2 < \dots < i_p.$ Here, we have defined $\gamma_{i}=b_i b_{i+g}$ for $1\le i  \leq g$, and require that 
\begin{equation}
m+2n+3p=3g-3.
\end{equation}\\

In \cite{Thaddeus} (after their equation (25), which is (\ref{twentyfive})
in our paper), Thaddeus writes: \emph{``Actually, the value of (25) is independent
of the choice of $i_j$. This follows from a diffeomorphism argument ... because
any permutation of the handles ... can be realized by a diffeomorphism.''}

We  define $$\gamma=2\sum_{k=1}^g \gamma_{k}. $$

The following fact is relevant:

\newcommand{\aaye}{a}
\newcommand{\beee}{f}
\newcommand{\psie}{b}
\newcommand{\psiee}{b}

\begin{prop} (\cite{Thaddeus}, Proposition 24):
The value of 
$$ \aaye^m \beee^n \prod_{i=1}^{2g} ( b_i)^{p_i} [\mathcal{M}_g]$$
is $0$ unless all the exponents $p_i$ appear in pairs $p_i = p_{i+g}$ and either
$p_i =0$ or $p_i = 1 $. \end{prop}

\begin{proof}

We start by showing that the quantity in question is zero unless:
\begin{enumerate} \item  $p_i \le 1$ for all $i$; 
\item 
$p_j = p_{j+g}  $ (and so these values must be either  $  0 $ or $1$).
\end{enumerate}

To prove the first claim, we observe that the $b_i$ are classes of odd degree 
so $(b_i)^2 = 0 $. 

To prove the second claim, we proceed as follows. 
Suppose $p_{i_0} + p_{i_0 + g} = 1$ (in other words one of the 
indices is $0$ and the other is $1$).
If we introduce an orientation-preserving
diffeomorphism $\psi$ on $\Sigma$, an isomorphism on cohomology will be induced, and 
so $\psi^* \aaye = \aaye$ and $\psi^* f = f$, 
but also 
$\psi^* b_{i_0}  =  - b_{i_0}   $ and
$\psi^* b_{i_0+g}  =  - b_{i_0+g}   $
while $\psi^* b_j = b_j $ for all the rest.
(A half twist of the surface around one loop will suffice.)
By naturality with respect to $\psi$, it follows  that
$$ \aaye^m \beee^n \prod_{i } (b_i )^{p_i} [\mathcal{M}_g]= -
 \aaye^m \beee^n \prod_{i } (b_i )^{p_i}[\mathcal{M}_g] = 0 .$$

\end{proof}

\begin{remark}
We note that the factors $2^g g!$ appear in the 
formulas for the intersection pairings. This follows from the fact that 
$\gamma/2$ is a sum of $g$ terms $b_i b_{i+g}$ as $i$ runs from $1$ to $g$.
As we proved above, when computing intersection pairings involving
powers of $\gamma$, we expand  the product $\gamma^n$ and the only terms
that contribute are those where each factor $b_i b_{i +g} $ appears 
at most once in the product. 
Any terms where this factor appears more than once 
contribute zero (as the degree of $b_j$ is odd and so it squares to $0$).

It follows that in the expansion of $\gamma^p$  there are
$2^p \left(\begin{array}{c}g \\ p\end{array}\right)$
terms, and 
all of these terms  contribute the same amount to the sum.

\end{remark}

\begin{theorem} \label{t3.1} (\cite{Thaddeus}, Proposition 26):

(a) Let $\eta_i: \mathcal{M}_{g-1} \to \mathcal{M}_g$ be the map induced by a map from $\Sigma_g$
to $\Sigma_{g-1}$ which collapses the $i$-th handle to a point. 

The image of $\mathcal{M}_{g-1}$ under $\eta_i$ is Poincar\'e dual to the class
$\gamma_i \in H^6(\mathcal{M}_g)$.

(b) Furthermore, we have 
\begin{equation} \label{eqn3}
a^m f^n \gamma^p [\mathcal{M}_g]=2g \,a^m f^n \gamma^{p-1} [\mathcal{M}_{g-1}].
\end{equation}
\end{theorem}

\begin{proof} 
First we prove part (a).
Recall that $M_g$ was defined in (\ref{defm}).
 Fix an element of $M_{g-1}$ of the form 

\begin{equation*}
(A_1, \dots, A_{i-1}, A_{i+1}, \dots, A_g, A_{g+1}, \dots ,A_{2g-2}),
\end{equation*}\\
and define a map $$\tau_i: SU(2) \to M_g$$ by 

\begin{align*}
\tau_i:   C &\mapsto (A_1, \dots, A_{i-1}, C,  A_{i+1}, \dots, A_g, A_{i+g-1}, C^{-1},A_{i+g},\dots,A_{2g-2}),
\end{align*}
which serves as a right inverse for the projection
map  $\pi_i: M_g \to SU(2)$ to the $i$-th component (for $1 \le i \le g$).
Likewise there is a right inverse to the projection map $\pi_{i+g}$. 
 Hence, if we pull back the cohomology classes $b_i$ by $\pi_i$ and also pull back $b_{i+g}$ by $\pi_{i+g}$, we will get indivisible cohomology classes in $M_g$ Poincar\'e dual to 
the homology cycles $\pi_i^{-1}(I)$ and $\pi_{i+g}^{-1}(I)$ of $M_g$, denoted respectively by 

\begin{equation*}
\chi_i, \chi_{i+g} \in H^3(M_g, \textbf{Z}).
\end{equation*}\\
If we denote the quotient map by $q: M_g \to \mathcal{M}_g$,
  we then claim that

\begin{align*}
\chi_i&=q^*(b_i) \\
\chi_{i+g}&=q^*(b_{i+g}) 
\end{align*}\\
for $\psi_i \in H^3(\mathcal{M}_g, $ 
\noin up to a sign (which 
turns out to be unimportant). This claim can be proved by assuming $\chi_i=q^*(\hat{b_i})$ for some indivisible $\hat{b_i}\in H^3(\mathcal{M}_g,\mathbb{Z})$. This $\hat{b_i}$ is indeed unique, and the only classes invariant under $f^*$ for all diffeomorphisms $f:\Sigma\to\Sigma$ are integer multiples of 
$b_i$. By the indivisibility of $\hat{b_i}$ and $b_i$, we can in fact conclude that $\hat{b_i}=\pm b_i$ which implies $\chi_i=\pm q^*(b_i)$.  

The proof of part (b)
 is immediate from part (a)  by definition of Poincar\'e duality.
\end{proof}

\begin{remark}
It follows from \eqref{eqn3} that all intersection pairings will be  known if we consider those
intersection pairings of the form
$$  \aaye^m \beee^n [\mathcal{M}_g].$$
Thaddeus gives a formula for these intersection pairings
using the Verlinde formula, a formula for the dimension of the 
space of holomorphic sections of  complex line bundles over $\mathcal{M}_g$, 
writing the Riemann-Roch formula in terms of the classes $a$ and $f$
(\cite{Thaddeus}, equation (29)):
\begin{equation} \label{verlinde}
 a^m f^n [\mathcal{M}_g] = (-1)^{g-1} \frac{m!}{(m-g+1)!} 2^{2g-2} 
(2^{m-g+1} -2) B_{m-g+1} 
\end{equation}
where $B_j$ is the $j$-th Bernoulli number.
\end{remark}

\begin{remark}
We observe finally that 
when $g = 2$ we have that Theorem \ref{t3.1} (b) 
(in other words, the second part of  Thaddeus,  Proposition 26) reduces to the fact that
 $\gamma_j[\mathcal{M}_{2} ]$
is  the number of points in $\mathcal{M}_{1}   $ (for $j = 1,2$)
and hence 
 $\gamma[\mathcal{M}_{2} ]$
is  $4$ times the number of points in $\mathcal{M}_{1}$.  In turn, $\mathcal M_1$ may be identified with 
$$ \{(x, y) \in SU(2)\times SU(2) \;|\; xyx^{-1}y^{-1} = - I\}/SU(2)$$
where the $SU(2)$ action is conjugation.  The quotient has a single orbit, represented by
$$ x  = \matr{i}{0}{0}{-i},\; y = \matr{0}{i}{i}{0}.$$It follows that 
$$  \gamma_j[\mathcal{M}_{2}  ] =1   ~~~~~~(j = 1,2)$$
and 
$$  \gamma[\mathcal{M}_{2}  ] = 4. $$
\end{remark}

\section{Poincar\'e Duals of Generators} \label{s:pd}
\vspace{2ex}

\noin Now that we  have in place  a complete set 
of generators for $H^*(\mathcal{M}_g)$, we would like to introduce systematic techniques for calculating the Poincar\'e duals of these generators. 

\subsection{Degree Two Generators} \hfill

\vspace{2ex}
The first generator of $H^*(\mathcal{M}_g)$ for which we would like to find its Poincar\'e dual is  $f\in H^2(\mathcal{M}_g)$. This generator is of particular interest: under an appropriate normalization, it represents the cohomology class of the symplectic form on $\mathcal{M}_g$. The details of this normalization can be found in \cite{JeffreyKirwan}, for example. The method of calculating the dual of this generator is via the zero section of a \textit{prequantum line bundle}. \\

\begin{definition} Let $(M,\omega)$ be a symplectic manifold. A \textbf{prequantum line bundle with connection} on $(M,w)$ is a complex line bundle $\mathcal{L} \to M$ equipped with a connection $\nabla$ for which the curvature $F_\nabla$ is equal to the symplectic form. i.e.\end{definition}

\begin{equation}
F_\nabla = \omega.
\end{equation} \\
\noin The Poincar\'e dual of $\omega$ may then be identified as the zero locus of a section into $\mathcal{L}$.
 It is worth noting here that all line bundles over $\mathcal{M}_g$ are simply powers of $\mathcal{L}$. This was proven in \cite{DrezNara}.

A construction of the prequantum line bundle can be found in \cite{ChernSimons}, and is expanded on in \cite{PQLB}. The authors of \cite{ChernSimons} construct this bundle in detail, and discuss sections of it. In \cite{PQLB} it is shown this line bundle has degree 1 when $\Sigma$ has
genus $1$,  and its connection to the symplectic form is elaborated on. 
Both these papers make the hypothesis that they are considering
flat connections on a trivial bundle over a 2-manifold. This is the 
case where the degree $d$ of the bundle is $0$, not the case where
the rank and degree are coprime.

We now describe an analogous situation where the Poincar\'e 
duals of degree two classes can be identified explicitly.
In \cite{weitsman}, Weitsman makes use of the fact that
some degree two cohomology classes are   Poincar\'e dual to the 
zero loci of  sections of
line bundles  to  explicitly  construct the intersection of
a number of such vanishing cycles and show that it is empty.
Specifically, Weitsman studies spaces $\mathcal{S}_g (t_1, \dots, t_N)$ 
of conjugacy classes of  flat $SU(2)$ 
connections over surfaces with $N$ boundary components
where the holonomy of the connections over each boundary
component is 
constrained to a constant value (the trace of the  holonomy around
the $m$-th boundary component
is $2 \cos (\pi t_m)$ for all $m = 1, \dots, N$).
The author  constructs circle bundles
$V^g_m (t_1, \dots, t_n) $ over $S_g(t_1, \dots, t_N)$  
and identifies their Poincar\'e duals (see 
\cite{weitsman}, Proposition 3.7). By geometric
arguments he is able to show that sufficiently high intersections
of sufficiently many of these Poincar\'e duals are empty.
By this means he gives a proof of an 
analogue of the conjecture of  Newstead \cite{Newstead1}
according to which $a^g = 0 $, where $a$ is the 
degree $4$ generator described above. 

This argument is adapted by Gamse and Weitsman \cite{GW} to show
vanishing of intersections of powers  of suitable cycles. 
This article generalizes the techniques of \cite{weitsman} to 
$SU(n)$. It gives a concise summary of the main argument of 
\cite{weitsman}. 
We specialize to 
 $SU(2)$, and restrict to a surface with one boundary
component where the holonomy around the boundary
component is $t \in SU(2)$, where $t$ is assumed not to 
be in the center $\{ \pm I\}$ .
Gamse and Weitsman define  a section which takes
a connection to the off-diagonal part of the $m$-th 
matrix $A_m$ in $(A_1, \dots, A_g, B_1, \dots, B_g)$. 
The Poincar\'e dual of 
the first Chern class of a line bundle 
they define is is the zero locus of this section.  This Poincar\'e 
dual 
consists of the connections for which the $m$-th matrix is in 
the maximal torus (in other words, it is diagonal). 
The authors find that if sufficiently many of these
zero loci have nonempty intersection, it follows that 
$$\prod_{i =1}^g [A_i, B_i] = I,$$ which is a contradiction 
since it is assumed that $$\prod_{i =1}^g [A_i, B_i] = t$$
where $t \ne I$ is an element of $SU(2)$.

\vspace{1ex}
\subsection{Degree Three Generators} \hfill

\vspace{2ex}
To construct the Poincar\'e duals of $b_1, \dots, b_{2g} \in H^3(\mathcal{M}_g)$, we first note that the $b_1, \dots, b_{2g}$ are integral generators of $H^3(\mathcal{M}_g, \textbf{Z})$. We find the Poincar\'e 
duals of these generators by working in the space of representations of $\pi_1(\Sigma)$ considered in section 2,
which is denoted by  $M_g$. 

In this setting, we can consider the operation of collapsing the $i^{th}$ handle of the orientable 2-manifold to a point. This induces an embedding

\begin{equation*}\eta_i : \mathcal{M}_{g-1} \hookrightarrow \mathcal{M}_g\end{equation*}

\noin of a lower genus moduli space into our original moduli space. More precisely, this embedding is the result of setting two of the generators of the fundamental group to the identity in $M_g$. If we define $\pi_i$ to be the restriction of $(SU(2))^{2g}$ to $M_g$ and projecting onto the $i^{th}$ component, then this becomes 

\begin{equation*}
\eta_i(\mathcal{M}_{g-1})=(\pi_i^{-1}(I) \cap \pi_{i+g}^{-1}(I))/SO(3).
\end{equation*} \\
We have dim$_\textbf{R}(\mathcal{M}_g)=6g-6$, and thus $\mathcal{M}_{g-1}$ will have real codimension $6$. We are now in a position to construct its Poincar\'e dual, which was the content of Theorem 
\ref{t3.1} (a) proved above. That result is due to Michael
Thaddeus  \cite{Thaddeus}.

In the proof of Theorem 
\ref{t3.1}(a) given above, it is shown that the 
Poincar\'e dual of the class $b_j$ 
pulls back under
the quotient map $q: M_g \to \mathcal{M}_g$ 
 to  the cycle  $\pi_j^{-1}(I)  \subset M_g$, 
where $\pi_j$ is the projection to the $j$-th copy of $SU(2)$.

\vspace{5pt}

\subsection{Degree Four Generators} \hfill

Finally, we wish to study the Poincar\'e dual of the degree four generator 
$$a\in H^4(\mathcal{M}_g).$$
 On our moduli space $\mathcal{M}_g,$ we may define characteristic  classes 
of the restriction  of  the  universal bundle $\mathbb{U}$ 
over $\mathcal{M}_g \times \Sigma$
to  $\mathcal{M}_g \times \{p\}  $, where  $p \in \Sigma$ is a point.

The class $a$ is 
the  Euler class of the 
restriction of  $\mathbb{U}$
to  $\mathcal{M}_g \times \{p\}$. 
The Euler class of a
 rank $2$ vector bundle is the same as the second 
Chern class.  Then, we have from \cite{BottandTu}: \\

\begin{prop} The Poincar\'e dual of the Euler class of a vector bundle $E \xrightarrow[]{\pi} M$ over an oriented manifold $M$ is the submanifold which is a zero section of E.\end{prop}

So the Poincar\'e dual of the 
degree four generator $a$ is the zero locus of 
a section of the bundle $\mathbb{U}$ restricted to 
$\mathcal{M}_g\times \{p\}$.

\section{Non-compact Analogue}

A natural question is how this extends to the non-compact setting, that is, to $SL(2,\mathbb C)$.  Here, the character variety itself is non-compact and is diffeomorphic to the moduli space of rank $2$ stable Higgs bundles via nonabelian Hodge theory.  The generators of this moduli space are identified in \cite{HT}.  In future work, we will consider the problem of explicitly calculating homological duals for these generators.\\

\textbf{Acknowledgements.}  The first and third named authors acknowledge NSERC for support through their respective Discovery Grants.  The second named author acknowledges the support of research assistantships from the University of Saskatchewan and the University of Toronto at Scarborough.  The authors are also grateful to an  anonymous referee for a careful
reading of the manuscript and many helpful suggestions.\\

\end{document}